# Asymptotic Intermediate Efficiency of the Chi-squared and Log-likelihood Ratio Goodness of Fit Tests

Sherzod M. Mirakhmedov

Institute of Mathematics. Academy of Sciences of Uzbekistan
100125, Tashkent, Mirzo-Ulugbek str., 81
e-mail: shmirakhmedov@yahoo.com

**Abstract**. The problem of testing the goodness of fit of an absolutely continuous distribution to a set of observations grouped into equal probability intervals, against to a family of sequences of alternatives approaching the hypothesis is considered. We focus our study to the asymptotic intermediate efficiency (AIE) due to Inglot (1999) of the chi-squared and log-likelihood ratio tests. It is shown that the AIE of these tests depend on the "rate of proximity" of the sequences of alternatives to the testing hypothesis. In particularly, it is determined a threshold family of alternatives, beyond of which the chi-squared test becomes inferior to the log-likelihood ratio test in terms of Inglot AIE. Asymptotical comparison of these tests through asymptotical behavior of the probabilities of type 1 errors, subject the tests have the same fixed asymptotic power is considered also. It is assumed that the number of groups increases together with the number of observations.

**Keywords.** Asymptotic efficiency, asymptotic slope, chi-square statistic, log-likelihood ratio statistic, multinomial distribution.

**MSC (2000):** 60F10, 62G10, 62G20

## 1. Introduction

The goodness-of-fit tests from groped data constitute a classical problem in statistical inference. Here the original problem is to test whether a sample has come from a given distribution. This problem is transformed into a problem of fit for a multinomial distribution by the method of grouping data: the support of given distribution is divided into $N$ mutually exclusive intervals and is observed the number of observations, say $\eta_k$, arisen in the $k$ th interval. Test statistics in this approach are indicate in different ways a divergence of observed variables $\eta_m$ from their expected values under the testing hypothesis. Most notable test statistics of this kind are Pearson's chi-squared ($\chi^2$) statistic:

$$\chi_N^2 = \sum_{m=1}^{N} \frac{(\eta_m - np_{0m})^2}{np_{0m}}$$

and the log-likelihood ratio (LR) statistic:

$$\Lambda_N = 2\sum_{m=1}^{N} \eta_m \ln \frac{\eta_m}{np_{0m}} \ ,$$

where $p_{0m}$ is the probability that the observation come from the $m$ th interval counted under the null hypothesis, $n$ is the sample size.

There is huge literature where interest and results have followed many aspects: the asymptotic distributional and statistical properties and recommendations in applications of these statistics in the case fixed $N$, the number of intervals, see Gvanceladze and Chibisov (1979), Moore (1986), Read and Cressie (1988), Cressie and Read (1989) and references within. However the assumption $N$ fixed is restrictive in several contexts. Indeed: Mann and Wald (1942) have obtained the relation $N \sim cn^{2/5}$, where $c > 0$ depends on the test size, concerning the optimal choice of the number of groups in chi-squared goodness of fit test; next, consider a problem of testing uniformity distribution on a fixed interval in which the interval partitioned into subintervals of equal length, if it is desired to achieve a specified expected frequency $v$ (say) for each subinterval, then $N$ subintervals are used for a sample size of $n$, where $N = N(n)$ is selected to make $n/N(n)$ close to $v$. For another examples associating with "big-data" problems see Rempata and Wesolowski (2016) and Pietrzak et al.( 2016).

In the present paper we study asymptotic properties in testing goodness-of-fit of $\chi_N^2$ and $\Lambda_N$ statistics in situation when $N = N(n) \to \infty$, as $n \to \infty$. Specifically, we consider the problem of testing the goodness of fit of an absolutely continuous distribution $F$ to a set of $n$ observations grouped into $N$ equal probability intervals. Through a probability integral transformation the original problem can be reduced to testing the null hypothesis of uniformity $H_0: F'(x) = f(x) = 1$, $0 \leq x \leq 1$. We test $H_0$ against a family of the sequences of alternatives:

$$H_{1n}: f_n(x) = 1 + \delta(n) l_n(x) \qquad (1.1)$$

where $\int_0^1 l_n(x) dx = 0$, $\|l_n\|_2 = 1$, $\sup_n \|l_n\|_\infty < \infty$ and $\delta(n) \to 0$ as $n \to \infty$.

Under the hypothesis $H_0$ the probabilities of intervals are $p_{01} = ... = p_{0N} = N^{-1}$, and hence, the $\chi_N^2$ and $\Lambda_N$ are special variants of the symmetric statistics:

$$S_N^h = \sum_{m=1}^N h(\eta_m), \qquad (1.2)$$

where $h$ is a function defined on non-negative axis. A test based on the statistic $S_N^h$ is called $h$-test for brevity. Asymptotical properties of $h$-tests essentially depend on behavior of $\delta(n)$, hence the alternatives (1.1) need some classification. Due to Holst (1972), Ivchenko and Medvedev (1978) and Mirakhmedov (1987) there is no power of the $h$-tests for the alternatives (1.1) with $\delta(n) = o\left((n\lambda_n)^{-1/4}\right)$, where $\lambda_n = n/N$. For the alternatives (1.1) with $\delta(n) = (n\lambda_n)^{-1/4}$ the chi-squared test is asymptotically most powerful within class of $h$-tests (satisfying some Lyapunov's kind condition). The rate $\delta(n) = (n\lambda_n)^{-1/4}$ keeps asymptotic power of the $h$-tests bounded away from the level and 1, and hence (1.1) in this case form a family of Pitman alternatives. Another family of "extreme" alternatives arises when one assumes that $\delta(n)$ remains fixed, i.e. alternatives does not approach the hypothesis; this case arises in Bahadur and Hodges-Lehman settings. Quine and Robinson (1985) have proved that when $\lambda_n \to \infty$ the chi-squared and LR tests have the same Pitman *asymptotical efficiency* (**AE**), but the chi-square test is inferior to the LR test in terms of the Bahadur AE. Further, these sequences of alternatives (1.1) with $\delta(n) \to 0$ and $\delta(n)(n\lambda_n)^{1/4} \to \infty$ provide family of intermediate alternatives, we shall denote it by symbol $\Im_{alt}$. According to our best knowledge there is only the paper by Ivchenko and Mirakhmedov (1995) who have considered the intermediate properties of $h$-tests in terms of asymptotical $\alpha$-slopes, see below Section 2, in testing hypothesis on the uniformity of the probabilities of a multinomial distribution. For the sparse multinomial model, i.e. $\lambda_n \to \lambda \in (0,\infty)$ and $0 < c_1 \leq Np_m \leq c_2 < \infty$, $m = 1,...,N$, they have proved that the chi-squared test is optimal in term of $\alpha$-slope within the class of $h$-tests for the subfamily of $\Im_{alt}$ such that $\delta(n) = O\left(\left((n\lambda_n)^{-1} \log n\right)^{1/4}\right)$, but for the subfamily of $\Im_{alt}$ satisfying the condition $\delta(n) n^{1/6} \log^{-1/3} N \to \infty$ the chi-square test is inferior to tests satisfying the Cramér condition, see Appendix, condition (A.1). Note that the $\Lambda_N$ statistic satisfies the Cramér condition, whereas the $\chi_N^2$ statistic does not.

Thus existing results does not cover the properties of the chi-square and LR tests for the family $\Im_{alt}$. Further, because the chi-square test being more (when $\lambda_n$ is bounded away from zero and infinity) or the same (if $\lambda_n \to \infty$) AE w.r.t. LR test for the Pitman family of alternatives (closest to hypothesis and distinguishable by the family of $h$-tests), and because it loses its leadership for the family of fixed alternatives, hence very interesting problem is to determine the family of alternatives at maximum "distance" from the hypothesis and so that the chi-squared test still retain "leadership". In this work among others we address those problems in terms of asymptotic relative efficiency, viz. a limiting value of the ratio of the sample sizes, developed by Inglot (1999) for the intermediate efficiency of tests goodness of fit; for the definitions see Section 2 and Appendix. Let $\Im_\gamma$ stands for the subfamily of $\Im_{alt}$ such that $\delta(n) = (n\lambda_n^2)^{-\gamma}$, $\gamma > 0$. In particularly, we show that for the families $\Im_\gamma$ with $\gamma > 1/8$ the chi-squared and LR tests are asymptotically same efficient, but for the families $\Im_\gamma$ with $0 < \gamma \leq 1/8$ the chi-squared test become inferior to the LR test. These complement the results of Quine and Robinson (1985), Mirakhmedov (1987) and Ivchenko and Mirakhmedov (1995).

The rest of the paper is organized as follows. The main results are presented in Section 2; in Section 3 the proofs are given; for the reader's convenience, the auxiliary Assertions are collected in Appendix. In what follows $c_j$ is a positive constant; all asymptotic statements are considered as $n \to \infty$; $\varsigma \sim F$ stands for "r.v $\varsigma$ has the distribution $F$"; $a_n \ll b_n$ stands for $a_n = o(b_n)$.

**2. The results**

We still use the notation of Sections 1. Remind, the symbol $\Im_{alt}$ stands for the family of alternatives (1.1), where $\delta(n) \to 0$ and $\delta(n)(n\lambda_n)^{1/4} \to \infty$, $\lambda_n = n/N$. It turned out that asymptotic properties of the chi-square test are different for the following subfamilies of the family of intermediate alternatives $\Im_{alt}$:

$\Im_o$, the subfamily such that $\delta(n) = o\left(\left(n \max(1, \lambda_n^2)\right)^{-1/6}\right)$,

$\Im_\gamma$, the subfamily such that $\delta(n) = (n\lambda_n^2)^{-\gamma}$, $1/8 < \gamma \leq 1/6$,

$\overline{\Im}_{1/8}$, the subfamily such that, $\delta(n) \geq (n\lambda_n^2)^{-1/8}$.

Before moving on to the specified AE results of the chi-square and LR tests we need in some notation on the general $h$-tests. Let $\alpha_n$ and $\beta_n$ stand for the size and power respectively of the $h$-test, $P_i$, $E_i$, $Var_i$ the probability, expectation and variance counted under $H_i$, $i = 0,1$; here and everywhere in the sequel $H_1$ means the family of alternatives under consideration, and $h$ is not a linear function. We suppose that large values of $S_N^h$ rejects the hypothesis and $E_1 S_N^h > E_0 S_N^h$.

Set $\xi \sim Poi(\lambda_n)$ and

$$\gamma = \lambda_n^{-1} \text{cov}(h(\xi), \xi), \quad \rho(S_N^h, \lambda_n) = corr\left(h(\xi) - \gamma\xi, \xi^2 - (2\lambda_n + 1)\xi\right). \tag{2.1}$$

Note that $\rho(\chi_N^2, \lambda_n) = 1$ and $\rho\left(S_N^h, \lambda_n\right) = corr_0(S_N^h, \chi_N^2)(1 + o(1))$, see Lemma 1 of Ivchenko and Mirakhmedov (1995).

Let $AIE\left(\chi_N^2, S_N^h\right)$ stands for the *asymptotic intermediate efficiency* (AIE) of $\chi_N^2$ test w.r.t. $S_N^h$ test due to Inglot (1999), see Definition A3 in Appendix.

**Theorem 2.1**. The following assertions are holding true:

a) For the families $\Im_o$ and arbitrary $\lambda_n$ one has

$$AIE\left(\chi_N^2, \Lambda_N\right) = \rho^{-2}\left(\Lambda_N, \lambda_n\right)(1 + o(1));$$

b) For the family $\Im_\gamma$, where for each $\gamma \in (1/8, 1/6]$ the $N$ such that

$$n^{(1-6\gamma)/(1-4\gamma)} \ll N \ll n^{3(1-4\gamma)/4(1-2\gamma)}, \tag{2.2}$$

one has

$$AIE(\chi_N^2, \Lambda_N) = 1;$$

c) For the family $\overline{\Im}_{1/8}$ if $\lambda_n \to \infty$ and $\delta(n) = o(\lambda_n^{-1/2})$ then

$$AIE(\chi_N^2, \Lambda_N) = 0.$$

**Theorem 2.2.** Let $\lambda_n \to \lambda \in (0, \infty)$. Then in the family $\Im_o$ for arbitrary $h$-test satisfying the Cramer condition (A.1), in particularly for LLR test, one has $AIE\left(\chi_N^2, S_N^h\right) = \rho^{-2}\left(S_N^h, \lambda\right)(1 + o(1)) > 1$.

Note that Theorem 2.2 deal with family $\Im_o$ only, this is because the region of the large deviation result for the chi-square statistics is limited by Assertion A2, see Appendix.

Ivchenko and Mirakhmedov (1995) have suggested as a measure of efficiency of $h$- test the asymptotic value of the slopes, which defined as follows.

Assume that $\beta_n \to \beta \in (0.1)$. The measure of performance of $h$-test is the asymptotic value of a $\alpha$-slope, viz.,

$$e_n^\alpha(S_N^h) = -\log P_0\{S_N^h \geq E_1 S_N^h\}. \tag{2.3}$$

Similarly, when $\alpha_n \to \alpha \in (0.1)$ the measure of performance of h-test is the asymptotic value of a $\beta$-slope, viz.,

$$e_n^\beta(S_N^h) = -\log P_1\{S_N^h \leq E_0 S_N^h\}. \tag{2.4}$$

We shall focus in the sequel on the $\alpha$-slope of $h$-tests. For a given family of alternatives the $\alpha$-asymptotic relative efficiency ($\alpha ARE$) of one test to another, subject they have the same asymptotic power, which is bounded away from 0 and 1, is defined as the ratio of their asymptotic $\alpha$-slopes. This corresponds to a concept of comparison of two tests via comparison of asymptotic behavior of the probabilities of the first type errors, when they have the same fixed asymptotic power; in our opinion such concept is reasonable and corresponds to the spirit of ARE of two tests.

The proof of Theorems 2.1 and 2.2 essentially uses the following results, which are of independent interest.

**Theorem 2.3**. The followings are hold:

a) In the family $\mathfrak{I}_o$ for arbitrary $\lambda_n$, and

b) In the family $\mathfrak{I}_\gamma$, where for each $\gamma \in (1/8, 1/6]$ the $N$ satisfies (2.2), one has

$$\frac{e_n^\alpha(\chi_N^2)}{n\lambda_n \delta^4(n)} = \frac{1}{4}(1+o(1)); \tag{2.5}$$

c) In the family $\overline{\mathfrak{I}}_{1/8}$ one has

$$\frac{e_n^\alpha(\chi_N^2)}{n\lambda_n \delta^4(n)} = o(1) \ .$$

**Theorem 2.4.** Let $\lambda_n \to \lambda \in (0,\infty]$. For the subfamily of $\mathfrak{I}_{alt}$ such that $\delta(n) = o(\lambda_n^{-1/2})$ one has

$$\frac{e_n^\alpha(\Lambda_N)}{n\lambda_n \delta^4(n)} = \frac{1}{4}\rho^2(\Lambda_N, \lambda)(1+o(1)) .$$

Note that $\rho(\Lambda_N, \lambda_n) = 1 - (6\lambda_n)^{-1}(1+o(1))$ if $\lambda_n \to \infty$. It follows from Theorem 2.3 and 2.4 that the $\alpha$ ARE of $\chi_N^2$ test w.r.t. $S_N^h$ test coincides with $AIE(\chi_N^2, S_N^h)$.

**Remark 2.1.** In Theorem 2.1 the case $\gamma > 1/6$ is in a) part. For $\gamma > 1/8$ the condition (2.2)

implies $N = o(\sqrt{n})$. In particularly, if $\gamma = 1/6$ then $N = o(n^{3/8})$, and if $\gamma = 1/8 + \varepsilon$, $\varepsilon \in (0, 1/24]$, then condition (2.2) implies $n^{(1-24\varepsilon)/2(1-8\varepsilon)} \ll N \ll n^{(3-24\varepsilon)/2(3-8\varepsilon)}$. So, $N$ is close to $\sqrt{n}$ for the $\gamma$ close to 1/8, and when $\gamma$ varies in the interval $(1/8, 1/6]$ the strips (2.2) cover interval $(1, o(\sqrt{n}))$.

**Remark 2.2.** Ivchenko and Mirakhmedov (1995) have considered $h$-tests in the problem of goodness of fit for a multinomial distribution $M(n, p_1, ..., p_N)$, namely the testing hypothesis $H_0$: $p_m = 1/N$, $m = 1, ..., N$, against sequences of alternatives $H_1: (p_1, ..., p_N) \neq (1/N, ..., 1/N)$, which approach $H_0$ so that $\varepsilon(N) = N^{-1} \sum_{m=1}^{N} (Np_m - 1)^2 \to 0$, $\sqrt{n\lambda_n} \varepsilon(N) \to \infty$. The above presented results can be reformulated and proved (by rewriting line by line of the proofs) for this goodness of fit problem by writing $\varepsilon(N)$ instead of $\delta^2(n)$.

**Remark 2.3.** Using asymptotical normality result of Mirakhmedov (2007) by standard algebra it can be shown that for Pitman alternatives the ARE of $h$-tests (satisfying Lyapunov kind condition, see below Lemma 3.1) w.r.t. chi-squared test is equal to $\rho^2(S_N^h, \lambda_n)$, see also Holst (1972) and Ivchenko and Medvedev (1978), where such result was obtained for sparse multinomial model, i.e. when $\lambda_n \to \lambda \in (0, \infty)$ and $0 < c_1 \leq Np_m \leq c_2 < \infty$, under some additional conditions for the function $h$. That is, for family of Pitman alternatives and the families of intermediate alternatives the ARE of $h$-test w.r.t. chi-square test is determined by the same characteristic $\rho(S_N^h, \lambda_n)$, the measure of correlation of the $h$-test statistic $S_N^h$ and $\chi_N^2$ statistic.

The above results allow making classification of the AE of the chi-squared test w.r.t. log-likelihood depending on the speed of proximity of alternatives to the hypothesis. Namely:

- Let $\lambda_n \to \infty$. Then the family of alternatives (1.1) with $\delta(n) = (n\lambda_n^2)^{-1/8}$ is a threshold in the AE of the chi-square test: for the alternatives which are at the "distance" $\delta(n) \geq (n\lambda_n^2)^{-1/8}$ from the hypothesis the AE of $\chi_N^2$ test is inferior to $\Lambda_N$ test, whereas for the alternatives which are at the "distance" $\delta(n) < (n\lambda_n^2)^{-1/8}$ (which include Pitman family of alternatives also) the $\chi_N^2$ and $\Lambda_N$ tests have the same AE. These complement the results of Queen and Robinson (1985).

- Let $\lambda_n \to \lambda$, $0 < \lambda < \infty$. Then for the family of alternatives with $\delta(n) = o(n^{-1/6})$ the $\chi_N^2$ test asymptotically more efficient w.r.t. all $h$-tests satisfying the Cramer condition (A.1), in particularly w.r.t. $\Lambda_N$ test. This complements the results of Holst (1972), Ivchenko and Medvedev (1978), Mirakhmedov (1987) and Ivchenko and Mirakhmedov (1995). However, in this case AE between the chi-squared and log-likelihood tests for the family of sequences of alternatives which are at "distance" $\delta(n) \geq n^{-1/6}$ is open problem yet. This is because there is no appropriate large deviation result for $\chi_N^2$ statistic in the area $[N^{1/6}, o(N^{1/2})]$.

**Remark 2.4.** Above considered intermediate approach based on the $\alpha$-slope is somewhere between Pitman and Bahadur approaches. Alike one can consider intermediate approach based on the $\beta$-slope, the case intermediate between Pitman and Hodges-Lehman settings, when the performance of a $h$-test based on second kind error asymptotic and is measured by the asymptotic value of $e_n^\beta(S_N^h)$, see (2.4). We state that under the Pitman families of alternatives and $\Im_{alt}$ the $e_n^\beta(S_N^h)$ and $e_n^\alpha(S_N^h)$ are asymptotically equivalent. In Inglot and Ledwina (2004, Remark 4) it was remarked that an approach based on the second kind error asymptotic is completely uninformative for Cramer-von Mises test; but it may rise a situation when under intermediate alternatives it is informative. The above statement confirms the last sentence.

### 3. Proofs

We still use the notation of Sections 1 and 2. First we consider the symmetric statistics $S_N^h$ (1.2). Remain that $\xi \sim Poi(\lambda_n)$. In addition to denotes (2.1) we set

$$g(\xi) = h(\xi) - Eh(\xi) - \gamma(\xi - \lambda_n) , \quad \sigma^2(h) = Var\, g(\xi) = Varh(\xi)\left(1 - corr^2\left(h(\xi), \xi\right)\right),$$

$$L_{3,N} = E|g(\xi)|^3 / \sigma^3(h)\sqrt{N} , \quad \Phi(u) = \int_{-\infty}^{u} e^{-t^2/2} dt .$$

**Lemma 3.1.** If $L_{3,N} \to 0$, then for any $f_n \in \Im_{alt}$ one has

$$P_i\left\{S_N^h < u\sqrt{Var_i S_N^h} + E_i S_N^h\right\} = \Phi(u) + o(1), \quad i = 0,1 ; \tag{3.1}$$

and

$$E_0 S_N^h = NEh(\xi)(1+o(1)), \quad Var_1 S_N^h = Var_0 S_N^h(1+o(1)) = N\sigma^2(h)(1+o(1)), \tag{3.2}$$

$$x_n(h) \stackrel{def}{=} \left(E_1 S_N^h - E_0 S_N^h\right) / \sqrt{Var_0 S_N^h} = \sqrt{\frac{n\lambda_n}{2}} \delta^2(n) \rho(S_N^h, \lambda_n)(1+o(1)) . \tag{3.3}$$

**Proof.** From Theorem of Mirakhmedov (2007) it follows that under the hypothesis $H_i$ the statistics $S_N^h$ has asymptotically normal distribution with expectation $A_{iN} = \sum_{m=1}^{N} E_i h(\xi_m)$ and variance $\sigma_{iN}^2 = \sum_{m=1}^{N} Var_i g(\xi_m)$, where $\xi_m \sim Poi(np_{im})$ and $i=0,1$. Hence by well-known theorem on convergence of moments (see Theorem 6.14 of Moran (1984))

$$E_i S_N^h = A_{iN} + o(1) \text{ and } Var_i S_N^h = \sigma_{iN}^2(1+o(1)) . \tag{3.4}$$

Therefore (3.1) follows. Note that in our situation the random vector of frequencies of the groups $(\eta_1, ..., \eta_N) \sim M(n, p_{i1}, ..., p_{iN})$, where $p_{0m} = 1/N$ under $H_0$, and under the alternatives (1.1)

$$p_{1m} = \int_{(m-1)/N}^{m/N} (1 + \delta(n)l_m(x)) dx = \frac{1}{N}(1 + \delta(n)l_{mN}) , \quad m = 1, ..., N .$$

Because of this applying Taylor expansion idea in the right-hand sides of equations in (3.4) after some computation we derive (3.2) and (3.3), see also Holst (1972), Ivchenko and Mirakhmedov (1995). Proof is completed.

**Remark 3.1**. Note that the condition $L_{3,N} \to 0$ is fulfilled for the chi-square and LR statistics iff $n\lambda_n \to \infty$.

**Lemma 3.2.** For every distribution $f_n \in \mathfrak{S}_{alt}$ and for arbitrary symmetric statistic $\hat{S}_N^h = (S_N^h - E_0 S_N^h)/\sqrt{Var_0 S_N^h}$ such that $L_{3,N} \to 0$ the condition (i) of Definition A.2 satisfies with $b_n(f_n) = \sqrt{\lambda_n/2}\,\delta^2(n)\rho(S_N^h, \lambda_n)$.

**Proof**. Due to Lemma 3.1 we have

$$P_1\left\{\left|\frac{\hat{S}_N^h}{x_n(h)} - 1\right| < \varepsilon\right\} = P_1\left\{\left|\frac{S_N^h - E_1 S_N^h}{\sqrt{Var_1 S_N^h}}\right| < \varepsilon x_n(h)\sqrt{\frac{Var_0 S_N^h}{Var_1 S_N^h}}\right\}$$

$$= 2\Phi\left(\varepsilon x_n(h)\sqrt{\frac{Var_0 S_N^h}{Var_1 S_N^h}}\right) - 1 + o(1) = 1 + o(1),$$

since (3.2) and the fact that $x_n(h) \to \infty$ for the family $\mathfrak{S}_{alt}$. Hence $\sqrt{n}b_n(f_n) = x_n(h)$. Lemma 3.2 follows, since (3.3).

We will prove Theorem 2.3 and 2.4 first. By (3.3) we have

$$e_n^\alpha(S_N^h) = -\log P_0\left\{\frac{S_N^h - NEh(\xi)}{\sigma(h)\sqrt{N}} \geq \sqrt{\frac{n\lambda_n}{2}}\delta^2(n)\rho(S_N^h, \lambda_n)(1 + o(1))\right\}. \qquad (3.5)$$

**Proof of Theorem 2.3**. For the chi-square statistic $h(u) = (u - \lambda_n)^2/\lambda_n$, $Eh(\xi) = 1$, $\sigma^2(h) = 2$ and $\rho(\chi_N^2, \lambda_n) = 1$, hence by (3.5)

$$e_n^\alpha(\chi_N^2) = -\log P_0\left\{\chi_N^2 > x_n\sqrt{2N} + N\right\}, \qquad (3.6)$$

where $x_n = \sqrt{n\lambda_n/2}\,\delta^2(n)(1 + o(1))$. Note that $x_n = o\left(\left(\sqrt{N}\min(1,\lambda_n)\right)^{1/3}\right)$ within the family $\mathfrak{S}_0$. Next, within family $\mathfrak{S}_\gamma$ we have: for each $\gamma \in (1/8, 1/6]$ the first and second relations of the condition (2.4) imply $x_n = o(\sqrt{N})$ and $N^{3/2}/\sqrt{n}x_n \to 0$ respectively. Therefore, proof of the cases a) and b) is concluded by applying Assertion A2 and Assertion A4 respectively in (3.6).

*Proof of part* c). Remark that in the family $\overline{\mathfrak{S}}_{1/8}$ we have $x_n \geq n^{1/4}/\sqrt{2}$, therefore $x_n = o(\sqrt{N})$ implies $\sqrt{n} = o(N)$. Hence Assertion A4 can't be used here. Instead we will prove the following

**Lemma 3.3**. Let $\sqrt{n} = o(N)$. In the families $\overline{\mathfrak{S}}_{1/8}$ one has $e_n^\alpha(\chi_N^2) = o(n\lambda_n\delta^4(n))$.

**Proof.** Let $v(n) = \left\lfloor \lambda_n + \sqrt{\lambda_n + n\lambda_n \delta^2(n)} \right\rfloor + 1$. Note that $E_0 \chi_N^2 = N(1+o(1))$, $Var_0 \chi_N^2 = 2N(1+o(1))$ and $x_n(h) = \sqrt{n\lambda_n/2}\,\delta^2(n)(1+o(1))$. Use these facts to get

$$P_0\{\chi_N^2 > E_1 \chi_N^2\} = P_0\{\chi_N^2 - E_0 \chi_N^2 > x_n(h)\sqrt{Var_0 \chi_N^2}\}$$

$$= P_0\left\{\sum_{m=1}^{N}\left((\eta_m - \lambda_n)^2 - \lambda_n\right) > n\lambda_n \delta^2(n)(1+o(1))\right\}$$

$$\geq P_0\left\{\sum_{m=2}^{N}\left((\eta_m - \lambda_n)^2 - \lambda_n\right) \geq 0 \Big/ \eta_1 = v(n)\right\} P_0\{\eta_1 = v(n)\}$$

$$= P_0\left\{\sum_{m=1}^{N-1}\left((\hat{\eta}_m - \lambda_n)^2 - \lambda_n\right) \geq 0\right\} P_0\{\eta_1 = v(n)\}. \tag{3.7}$$

Here $\hat{\eta}_m \sim Bi(n - v(n), (N-1)^{-1})$. Put $\hat{\lambda}_n = (n - v(n))/(N-1)$. It is easy to see that $v(n)/n = (N^{-1} + \delta(n)N^{-1/2})(1+o(1))$ and $\hat{\lambda}_n = \lambda_n(1 + O(N^{-1} + \delta(n)N^{-1/2}))$. We have

$$P_0\left\{\sum_{m=1}^{N-1}\left((\hat{\eta}_m - \lambda_n)^2 - \lambda_n\right) \geq 0\right\} \geq P_0\left\{\sum_{m=1}^{N-1}(\hat{\eta}_m - \hat{\lambda}_n)^2 \geq (N-1)\lambda_n\right\}$$

$$= P_0\left\{\sum_{m=1}^{N-1}\left((\hat{\eta}_m - \hat{\lambda}_n)^2 - \hat{\lambda}_n\right) \geq (v(n) - \lambda_n)\right\}$$

$$= P_0\left\{\frac{\sum_{m=1}^{N-1}\left((\hat{\eta}_m - \hat{\lambda}_n)^2 - \hat{\lambda}_n\right)}{\sqrt{2(n - v(n))^2/(N-1)}} \geq \delta(n) + o(1)\right\} \geq c > 0, \tag{3.8}$$

because $(n - v(n))\hat{\lambda}_n = n\lambda_n(1+o(1)) \to \infty$, and hence the CLT for the statistic $\sum_{m=1}^{N-1}(\hat{\eta}_m - \hat{\lambda}_n)^2$ is enable to use, see Mirakhmedov (1992, Corollary 3).

Set $g(x, p) = x\log(x/p) + (1-x)\log((1-x)/(1-p))$, $x \in (0,1)$ and $p \in (0,1)$. Let $\varsigma \sim Bi(k, p)$. Due to Lemma 1 of Quine and Robinson (1985): for an integer $kx$

$$P\{\varsigma = kx\} \geq 0.8(2\pi kx(1-x))^{-1/2} \exp\{-kg(x, p)\}. \tag{3.9}$$

Note that under $H_0$ the r.v. $\eta_1 \sim B(n, N^{-1})$, therefore applying (3.9) we obtain

$$P_0\{\eta_1 = v(n)\}$$

$$\geq c\left(v(n)(1 - v(n)n^{-1})\right)^{-1/2} \exp\left\{-v(n)\log(\lambda_n^{-1} v(n)) - n(1 - n^{-1}v(n))\log\frac{1 - v(n)n^{-1}}{1 - N^{-1}}\right\}$$

$$\geq c\left(v(n)\right)^{-1/2} \exp\left\{-v(n)\log(\lambda_n^{-1} v(n))\right\}.$$

Hence

$$-\frac{\log P_0\{\eta_1 = v(n)\}}{n\lambda_n \delta^4(n)} \leq c \frac{\log v(n) + v(n)\log(\lambda_n^{-1} v(n))}{n\lambda_n \delta^4(n)} \leq c \frac{\lambda_n + \delta(n)\sqrt{n\lambda_n}}{n\lambda_n \delta^4(n)} \log \frac{v(n)}{\lambda_n}$$

$$\leq c\left(\frac{1}{n\delta^4(n)} + \frac{1}{\delta^3(n)\sqrt{n\lambda_n}}\right) \log \frac{v(n)}{\lambda_n} \leq c\left[\left(\frac{n}{N^2}\right)^{1/2} + \left(\frac{n}{N^2}\right)^{1/8}\right] \log\left(\frac{N^2}{n}\right) = o(1), \qquad (3.10)$$

since $\delta(n) \geq (N^2/n^3)^{1/8}$ within the families $\bar{\mathfrak{F}}_{1/8}$, and $n = o(N^2)$. Lemma 3.3 follows from (3.6), (3.7), (3.8) and (3.10).

For $\sqrt{n} = o(N)$ part c) follows from Lemma 3.3. Let now $N = o(\sqrt{n})$. Lemma 3.2 together with Assertion A4 allow to conclude that $e_n^\alpha(\chi_N^2)$ asymptotically coincides with intermediate slope, see Definition A2; here the constant $c$ and the sequence $\tau_n$ (from Definition A2) are such that $\tau_n c = e_n^\alpha(\chi_N^2)/n\lambda_n \delta^4(n)/4(1+o(1))$, since (3.3). On the other hand, when testing uniformity against family $\mathfrak{F}_{alt}$ the Neyman-Pearson (NP) test can be applied. As it was shown by Inglot and Ledvina (1996, Sec 5) the intermediate slope of NP test has the form $n\delta^2(n)/2$. Hence, if $\tau_n = 1$ under the family $\bar{\mathfrak{F}}_{1/8}$ then $e_n^\alpha(\chi_N^2)/e_n^\alpha(NP) = 2c\lambda_n \delta^2(n) \geq 2c(n/N^2)^{1/4} \to \infty$, but this is impossible, because the ratio $e_n^\alpha(\chi_N^2)/e_n^\alpha(NP)$ can't be greater than 1. Hence $\tau_n = o(1)$. Proof of Theorem 2.1 is completed.

**Proof of Theorem 2.4.** By (3.2), (3.3) and (3.5) $e_n^\alpha(\Lambda_N) = -\log P_0\left\{\Lambda_N > x_n \sqrt{Var_0 \Lambda_N} + E_0 \Lambda_N\right\}$, where $x_n = \sqrt{n\lambda_n/2}\,\delta^2(n)\rho(\Lambda_N, \lambda_n)(1+o(1))$. Note that $x_n = o(N^{1/2})$ iff $\delta(n) = o(\lambda_n^{-1/2})$. Proof of Theorem 2.4 is concluded by applying Assertions A3 and A5.

**Proof of Theorem 2.1 and 2.2.** Let $c$ be an arbitrary constant and
$t_n^h(f_n) = \sqrt{n\lambda_n/2}\,\delta^2(n)\rho(S_N^h, \lambda_n) + c$, $f_n \in \mathfrak{F}_{alt}$. The test statistic is $\hat{S}_N^h = (S_N^h - E_0 S_N^h)/\sqrt{Var_0 S_N^h}$. For the power of the $h$-test with critical region

$$\{\hat{S}_N^h > t_n^h(f_n)\} \qquad (3.11)$$

using Lemma 3.1 we obtain

$$P_1\{\hat{S}_N^h > t_n^h(f_n)\} = P_1\left\{\frac{S_N^h - E_1 S_N^h}{\sqrt{Var_1 S_N^h}} > t_n^h(f_n) - \frac{E_1 S_N^h - E_0 S_N^h}{\sqrt{Var_0 S_N^h}}\sqrt{\frac{Var_0 S_N^h}{Var_1 S_N^h}}\right\} = \Phi(-c) + o(1), \qquad (3.12)$$

since (3.2) and (3.3). Hence asymptotic power of this $h$-test is bounded away from zero and 1. The related significance level of (3.11) is defined via

$$\alpha_n(S_N^h) = P_0\{\hat{S}_N^h > t_n^h(f_n)\} = P_0\{\hat{S}_N^h > \sqrt{n\lambda_n/2}\,\delta^2(n)\rho(S_N^h, \lambda_n)(1+o(1))\}. \qquad (3.13)$$

Further proof is concluded by using (3.12), (3.13) facts in checking the conditions of Assertion A.6 with $V_n^{(2)}(f_n) = \hat{\chi}_N^2$ and $V_n^{(1)}(f_n) = \hat{\Lambda}_n$, the standardized version of the $\chi_N^2$ and $\Lambda_N$ statistics respectively, which are considered as the test statistics. Let $\mathbf{IS}(S_N^h)$ stands for the intermediate slope of the statistic $S_N^h$, see Definition A2. The condition (i) of the Definition A2 follows from Lemma 3.2. Next, from large deviations results, given in Assertions A1-A5, and Theorem 2.3 it follows that for $\hat{S}_n^f$ statistic and its special cases $\hat{\chi}_N^2$ and $\hat{\Lambda}_n$ the condition (ii) of Definition A2 is fulfilled for the appropriate families of alternatives. Therefore we obtain the following facts:

(a) The $\hat{\Lambda}_n$ is $(\Im_{alt}, q_n^{(1)}, 1)$-regular for the $\lambda_n \to \lambda \in (0, \infty]$, where $q_n^{(1)} = \lambda_n^{-1/2}$, and $\mathbf{IS}(\hat{\Lambda}_N) = 4^{-1} \lambda_n \delta^4(n) \rho^2(\Lambda_N, \lambda_n)$, remind that $\rho^2(\Lambda_N, \lambda_n) \to 1$ if $\lambda_n \to \infty$.

(b) The $\hat{\chi}_N^2$ statistic is $(\Im, q_n^{(2)}, 1)$-regular, where: (i) $\Im = \Im_o$, $q_n^{(2)} = n^{-1/3} \min(\lambda_n^{-1/6}, \lambda_n^{1/2})$ for arbitrary $\lambda_n$, and (ii) $\Im = \Im_\gamma$, $q_n^{(2)} = \lambda_n^{-1/2}$ for each $\gamma \in (1/8, 1/6]$ and $N$ from the strip (2.4); for these both cases $\mathbf{IS}(\hat{\chi}_N^2) = 4^{-1} \lambda_n \delta^4(n)$. Also, the $\hat{\chi}_N^2$ statistic is $(\overline{\Im}_{1/8}, \lambda_n^{-1/2}, \tau_n)$-regular with $\tau_n = o(1)$ and $\mathbf{IS}(\hat{\chi}_N^2) = \tau_n \lambda_n \delta^4(n)$.

(c) If $\lambda_n \to \lambda \in (0, \infty)$ then according to Assertion A1 the test statistic $S_N^h$ satisfying the Cramer condition (A.1) is $(\Im_{alt}, 1, 1)$-regular and $\mathbf{IS}(S_N^h) = 2^{-1} \sqrt{\lambda_n} \delta^2(n) \rho(S_N^h, \lambda_n)$.

Further, observe that $\Im_{alt}$ is renumerable, $\log \alpha_n(\chi_N^2) = -n \lambda_n \delta^4(n)/4(1+o(1))$ and $e = 1$ for the family $\Im_o$ and $\Im_\gamma$, $\gamma \in (1/8, 1/6]$, whereas $\log \alpha_n(\chi_N^2) = o(n \lambda_n \delta^4(n))$ and $e = 0$ for the family $\overline{\Im}_{1/8}$, since (3.13) and Theorem 2.3 (because $\log \alpha_n(S_N^h) = e_n^\alpha(S_N^h)$ by (3.3)). Due to these, (3.12) and above (a),(b) and (c) facts we see that all conditions of Assertion A6 are satisfied. Theorems 2.1 and 2.2 follow.

In order to confirm the statement of Remark 3.4 we note that

$$e_n^\beta(S_N^h) = -\log P_1 \left\{ \frac{S_N^h - E_1 S_N^h}{\sqrt{Var_1 S_N^h}} \leq -\sqrt{\frac{n\lambda_n}{2}} \delta^2(n) \rho(S_N^h, \lambda_n)(1+o(1)) \right\},$$

since (3.2) and (3.3). This implies that $\beta$-slope also is determined by the large deviation probabilities for the left tails of $S_N^h$ under alternatives $\Im_{alt}$. Such kind results are given in Mirakhmedov (2016) for the $\chi_N^2$ and $\Lambda_N$ statistics and by Ivchenko and Mirakhmedov (1995) for the some class of symmetric statistics.

**Appendix**. We still use the notation of the previous sections.

**Assertion A1**. Let the function $h$ be not linear, $\lambda_n \to \lambda \in (0, \infty)$, $Np_{\max} < c_1$, for some $c_1 > 0$, and

$$\max_{1 \leq m \leq N} E \exp\{H |g(\xi_m)|\} < c_2, \text{ for some } H > 0 \text{ and } c_2 > 0, \xi_m \sim Poi(np_m). \quad (A.1)$$

Then for $x_n \to \infty$, $0 \leq x_n = o(N^{1/2})$ it holds

$$\log P\left\{S_N^h > x_n \sqrt{VarS_N^h} + ES_N^h\right\} = -\frac{1}{2}x_n^2 + O\left(\log x_n + \frac{x_n^3}{\sqrt{N}}\right).$$

Assertion A1 follows from Theorem 2 of Ivchenko and Mirakhmedov (1995) and the fact that

$$1 - \Phi(x_n) = (x_n\sqrt{2\pi})^{-1} \exp\{-x_n^2/2\}(1+o(1)), \ x_n \to \infty. \tag{A.2}$$

Note that for the statistics $\Lambda_N$ Cramer condition (A.1) is satisfied but for the $\chi_N^2$ does not.

**Assertion A2.** Let $p_m = N^{-1}(1 + \delta(n)d_{m,n})$, $m = 1,2,...,N$, where $\delta(n) \to 0$ and

$$\sum_{m=1}^{N} d_{m,n} = 0, \ \frac{1}{N}\sum_{m=1}^{N} d_{m,n}^2 = d^2 < \infty,$$

Then for arbitrary $\lambda_n$ and $x_n$ such that $0 \leq x_n = o\left(\left(\sqrt{N}\min(1,\lambda_n^2)\right)^{1/3}\right)$ one has

$$P\left\{\chi_N^2 > x_n\sqrt{2N} + N\right\} = -\frac{1}{2}x_n^2 + O(\log x_n),$$

**Assertion A3.** Let $p_m = N^{-1}(1+o(1))$, $m = 1,2,...,N$, If $\lambda_n \to \infty$, $x_n \to \infty$, $x_n = o(N^{1/6})$ then

$$\log P\left\{\Lambda_N > x_n\sqrt{2N} + N\right\} = -\frac{1}{2}x_n^2 + O(\log x_n).$$

Assertion A2 and Assertion A3 follow from Corollary 2.2 and Corollary 2.3 of Mirakhmedov (2016) respectively, by applying (A.2).

**Assertion A4.** Let $N = o(\sqrt{n})$, $Np_{\min} > c$. If $x_n \to \infty$, $x_n = o(\sqrt{N})$ and $N^{-3/2}n^{1/2}x_n \to \infty$ then

$$\log P\left\{\chi_N^2 > x_n\sqrt{2N} + N\right\} = -\frac{1}{2}x_n^2 + O\left(\frac{x_n^3}{\sqrt{N}} + \log N + \frac{x_n N^{3/2}}{\sqrt{n}}\right).$$

**Assertion A5.** Let $\lambda_n \to \infty$ and $Np_{\min} > c$ for some $c > 0$. If $x_n \to \infty$ and $x_n = o(\sqrt{N})$ then

$$\log P\left\{\Lambda_N > x_n\sqrt{2N} + N\right\} = -\frac{1}{2}x_n^2 + O\left(\frac{x_n^3}{\sqrt{N}} + \log N + \frac{N^{3/2}}{\sqrt{n}}\right).$$

Assertions A4 and A5 are respectively Eq. (2.17) and Eq. (2.13) of Kallenberg (1985).

The notion of Kallenberg's intermediate efficiency, Kallenberg (1983), for the case when dealing with nonparametric sequences of alternatives has been developed by Inglot (1999). We need a slightly more general definition than his one. Below the notation of Inglot (1999) in the variant adapted to the problem of testing uniformity against sequences of alternatives (1.1) are presented.

So, we consider testing the hypothesis $H_0: f(x) = 1$, $0 < x < 1$. In what follows symbol $F^*$ denotes the family of all sequences of alternative distributions whose densities have the form (1.1), where $l_n(x)$ as in (1.1), but

$$\delta(n) \to 0 \text{ and } \sqrt{n}\delta(n) \to \infty, \tag{A.3}$$

**Definition A1.** We say that a family $F$ of alternatives (1.1) satisfying (A.3) is renumarable if for every $f_n \in F$ and all sequences of positive integers $n_j$ and $k_j$, such that $n_j \to \infty$, $n_j = O(k_j)$ the sequence $f'_n = f_n$ if $n \neq k_j$ else $f'_n = f_{n_j}$ if $n = k_j$, $j = 1, 2, ...$ also belong to $F$.

Note that the families $F^*$ and $\Im_{alt}$ are renumerable, whereas subfamilies $\Im_o, \Im_\gamma$ and $\bar{\Im}_{1/8}$, see Section 3, are not.

Let $q = (q_n)$ and $\tau = (\tau_n)$ be sequences of positive numbers such that

$$\lim q_n = 0, \lim n q_n^2 = \infty, \text{ or } q_n = 1 \text{ for all } n, \tag{A.4}$$

$$\lim \tau_n = 0, \text{ or } \tau_n = 1 \text{ for all } n. \tag{A.5}$$

Also, let $F$ be a subfamily of the family $F^*$, and $V_n(f_n)$, $f_n \in F$, be a test statistic rejecting the hypothesis for large values of $V_n(f_n)$.

**Definition A2.** The test statistic $V_n(f_n)$ is called $(F, q, \tau)$-regular for testing $H_0$ against $H_1$ if $V_n(f_n)$ is defined for every $f_n \in F$ and the following two conditions hold:

(i) there exists a positive function $b_n(f_n)$ such that for every $f_n \in F$ and arbitrary small $\varepsilon > 0$

$$\lim P_1 \left\{ \left| \frac{V_n(f_n)}{\sqrt{n} b_n(f_n)} - 1 \right| \leq \varepsilon \right\} = 1,$$

(ii) there exists a constant $c > 0$ such that for every sequence $x_n$, $x_n = o(q_n)$ and $\tau_n n x_n^2 \to \infty$

$$\lim \frac{1}{\tau_n n x_n^2} \log P_0 \left\{ V_n(f_n) \geq \sqrt{n} x_n \right\} = -c.$$

The function $c\tau_n [b_n(f_n)]^2$, $f_n \in F$, is called the intermediate slope of $V_n(f_n)$.

For two sequences of test statistics $V_n^{(1)}(f_n)$ and $V_n^{(2)}(f_n)$ with right critical regions let $\sqrt{n} t_{n,\alpha}^{(1)}$ and $\sqrt{n} t_{n,\alpha}^{(2)}$ be corresponding critical values of the same level $\alpha$, i.e. $P_0 \{V_n^{(i)} \geq \sqrt{n} t_{n,\alpha}^{(i)}\} \leq \alpha$, and $P_0 \{V_n^{(i)} \geq d\} > \alpha$ for every $d < \sqrt{n} t_{n,\alpha}^{(i)}$, $i = 1, 2$. Let $V_n^{(j)}(f_n)$ be $(F^{(j)}, q^{(j)}, \tau^{(j)})$-regular, where $F^{(j)} \subseteq F^*$, the sequences $q^{(j)} = (q_{jn})$, $\tau^{(j)} = (\tau_{jn})$ satisfy (A.4) and (A.5) respectively, $j = 1, 2$. Let $\alpha_n$ be a sequence of levels such that

$$\lim \alpha_n = \lim (\tau'_n n)^{-1} \alpha_n = 0, \text{ with } \tau'_n = \min(\tau_{1n}, \tau_{2n}), \tag{A.6}$$

and the power of $V_n^{(2)}(f_n)$ test of the level $\alpha_n$ is bounded away from 0 and 1:

$$0 < \liminf P_1 \{V_n^{(2)}(f_n) \geq \sqrt{n} t_{n,\alpha_n}^{(2)}\} \leq \limsup P_1 \{V_n^{(2)}(f_n) \geq \sqrt{n} t_{n,\alpha_n}^{(2)}\} < 1. \tag{A.7}$$

Put

$$n_{V^{(2)}, V^{(1)}} = n_{V^{(2)}, V^{(1)}}(n, f_n, \alpha_n) = \inf \left( m : P_1 \{V_{m+k}^{(1)}(f_n) \geq \sqrt{m+k} \, t_{m+k,\alpha_n}^{(1)}\} \right.$$

$$\geq P_1\left\{V_n^{(2)}(f_n) \geq \sqrt{n}\, t_{n,\alpha_n}^{(2)}\right\} \text{ for all } k \geq 0\Big).$$

**Definition A3**. Let a family $F$ satisfies (A.3) and there exists $\alpha_n$ satisfying (A.6) and (A.7). If there exists the limit

$$\lim \frac{n_{V^{(2)},V^{(1)}}}{n} = e_{V^{(2)},V^{(1)}} \in [0,\infty],$$

which does not depend on the particular choice of $\alpha_n$, then $e_{V^{(2)}V^{(1)}}$ is called asymptotic intermediate efficiency of $V_n^{(2)}(f_n)$ with respect to $V_n^{(1)}(f_n)$, shortly $\text{AIE}\left(V_n^{(2)}(f_n), V_n^{(1)}(f_n)\right)$.

**Assertion A6**. Let (i) aforementioned statistic $V_n^{(j)}(f_n)$ be $(F^{(j)}, q^{(j)}, \tau^{(j)})$-regular with corresponding functions $b_n^{(j)}(f_n)$ and constants $c^{(j)}$, $j = 1, 2$; (ii) $F^{(1)} \supseteq F^{(2)}$ and $q_{2n} \leq q_{1n}\sqrt{\tau_{1n}}$, for $n$ sufficiently large; (iii) $F^{(1)}$ be renumerable, $q_{1n}^{-1}$ and $\sqrt{\tau_n n}\, q_{1n}$ are non-decreasing; (iv) for every $f_n \in F^{(2)}$ there exists the limit

$$\lim \frac{c^{(2)} \tau_{2n} \left[b_n^{(2)}(f_n)\right]^2}{c^{(1)} \tau_{1n} \left[b_n^{(1)}(f_n)\right]^2} = e \in [0,\infty];$$

(v) for each $f_n \in F^{(2)}$ there exists $\alpha_n$ satisfying (A.6) and (A.7) and such that $\log \alpha_n = o(nq_{2n}^2)$. Then $e_{V^{(2)},V^{(1)}} = e$.

The proof of Assertion A6 consists of rewriting line by line the proof of Theorem 2.7 of Inglot (1999) with quite obvious changes, and is here omitted.

**References.**

1. Cressie N.A.C. and Read T.R. C. (1989). Pearson's $X^2$ and the Log-likelihood Ratio Statistic $G^2$: A Comparative Review. International Statistical Review, v.57, p.19-43.
2. Gvanceladze L.G. and Chibisov D.M. (1979). On tests of fit based on grouped data. In Contribution to Statistics, J.Hajek Mamorial Valume. J.Jurechkova. ed. 79-89. Academia, Prague.
3. Holst L. (1972). Asymptotic normality and efficiency for certain goodness-of-fit tests. Biometrica, **59**, p.137-145.
4. Inglot T. and Ledwina T.,(1996). Asymptotic optimality of data-driven Neyman's tests for uniformity. Ann. Statist. **24**, p.1982-2019.
5. Inglot T., (1999). Generalized intermediate efficiency of goodness-of-fit tests. Math. Methods Statist; **8**, p.487-509.
6. Inglot T. and Ledwina T.,(2004). On consistent minimax distinguishability and intermediate efficiency of Cramer-von Mises test. J. Stat. Planning Infer., **124**, p. 453-474.

7. Ivchenko G.I. and Medvedev Y.I.(1978). Decomposable statistics and verifying of tests. Small sample case. Theory of Probability Appl., **23** , 796-806.

8. Ivchenko G.I. and Mirakhmedov Sh.A., (1995). Large deviations and intermediate efficiency of the decomposable statistics in multinomial scheme. Math. Methods in Statist., **4**, p.294-311.

9. Kallenberg W.C.M.,(1983). Intermediate efficiency, theory and examples. Ann. Statist, **11**, p. 170-182.

10. Kallenberg W.C.M.,(1985). On moderate and large deviations in multinomial distributions. Ann. Statist. **13**, 1554-1580.

11. Mann H.B. and Wald, A. (1942). On the choice of the number of intervals in the application of the chi-square test. Ann. Math. Statist. **13**, 306-317.

12. Mirakhmedov S. A., (1987). Approximation of the distribution of multi-dimensional randomized divisible statistics by normal distribution.(Multinomial case). Theory Probabl. Appl., 32, p.696-707.

13. Mirakhmedov S .A.(1992), Randomized decomposable statistics in the scheme of independent allocating particles into boxes. Discrete Math. Appl. **2**, p.91-108. DOI: 10.1515/dma.1992.2.1.91, October 2009

14. Mirakhmedov S.M[1]. (2007). Asymptotic normality associated with generalized occupancy problem. Statistics & Probability Letters, v.77, p.1549-1558

15. Mirakhmedov S.M. (2016), The Probabilities of Large Deviations for Chi-square and Log likelihood Ratio Statistics. http://arxiv.org/abs/1606.00250.

16. Moore D.S (1986). Tests of chi-squared type in goodness of fit techniques. Ed. R.B. D'Agostino and V.A. Stephens, pp 63-95. New York: Marcel Dekker.

17. Moran P.A.P. (1984) An introduction to probability theory. Oxford University Press, Oxford, p 550

18. Quine M.P. and Robinson J. (1985). Efficiencies of chi-square and likelihood ratio goodness-of-fit tests, The Annals of Statistics. **13**, p.727–742.

19. Pietrzak, M., Rempała, G.A., Seweryn, M., Wesołowski, J.( 2016). Limit theorems for empirical Rènyi entropy and divergence with applications to molecular diversity analysis. TEST, 1–20.

20. Read T. R. C. and Cressie N. A. C. (1988). Goodness-of-Fit Statistics for Discrete Multivariate Data. Springer-Verlag, New York.

21. Rempata and Wesolowski (2016). Double asymptotics for the chi-square statistic. Statis. Probabl Letters, 119,p.317-325.

---

[1] Mirakhmedov S.M. is former Mirakhmedov Sh.A.